\documentclass[12pt]{article}
\usepackage{amsfonts}

\def\hybrid{\topmargin 0pt      \oddsidemargin 0pt
        \headheight 0pt \headsep 0pt
        \voffset=-0.5cm
        \textwidth 6.5in        
        \textheight 9in         
        \marginparwidth 0.0in
        \parskip 5pt plus 1pt   \jot = 1.5ex}
\catcode`\@=11
\def\marginnote#1{}

\newcount\hour
\newcount\minute
\newtoks\amorpm
\hour=\time\divide\hour by60
\minute=\time{\multiply\hour by60 \global\advance\minute by-\hour}
\edef\standardtime{{\ifnum\hour<12 \global\amorpm={am}%
        \else\global\amorpm={pm}\advance\hour by-12 \fi
        \ifnum\hour=0 \hour=12 \fi
        \number\hour:\ifnum\minute<10 0\fi\number\minute\the\amorpm}}
\edef\militarytime{\number\hour:\ifnum\minute<10 0\fi\number\minute}

\def\draftlabel#1{{\@bsphack\if@filesw {\let\thepage\relax
   \xdef\@gtempa{\write\@auxout{\string
      \newlabel{#1}{{\@currentlabel}{\thepage}}}}}\@gtempa
   \if@nobreak \ifvmode\nobreak\fi\fi\fi\@esphack}
        \gdef\@eqnlabel{#1}}
\def\@eqnlabel{}
\def\@vacuum{}
\def\draftmarginnote#1{\marginpar{\raggedright\scriptsize\tt#1}}

\def\draftlabel#1{{\@bsphack\if@filesw {\let\thepage\relax
   \xdef\@gtempa{\write\@auxout{\string
      \newlabel{#1}{{\@currentlabel}{\thepage}}}}}\@gtempa
   \if@nobreak \ifvmode\nobreak\fi\fi\fi\@esphack}
        \gdef\@eqnlabel{#1}}
\def\@eqnlabel{}
\def\@vacuum{}
\def\draftmarginnote#1{\marginpar{\raggedright\scriptsize\tt#1}}

\def\draft{\oddsidemargin -.5truein
        \def\@oddfoot{\sl preliminary draft \hfil
        \rm\thepage\hfil\sl\today\quad\militarytime}
        \let\@evenfoot\@oddfoot \overfullrule 3pt
        \let\label=\draftlabel
        \let\marginnote=\draftmarginnote
   \def\@eqnnum{(\theequation)\rlap{\kern\marginparsep\tt\@eqnlabel}%
\global\let\@eqnlabel\@vacuum}  }

\def\numberbysection{\@addtoreset{equation}{section}
        \def\theequation{\thesection.\arabic{equation}}}

\def\underline#1{\relax\ifmmode\@@underline#1\else
        $\@@underline{\hbox{#1}}$\relax\fi}

\def\titlepage{\@restonecolfalse\if@twocolumn\@restonecoltrue\onecolumn
     \else \newpage \fi \thispagestyle{empty}\c@page\z@
        \def\thefootnote{\fnsymbol{footnote}} }

\def\endtitlepage{\if@restonecol\twocolumn \else  \fi
        \def\thefootnote{\arabic{footnote}}
        \setcounter{footnote}{0}}  
\relax

\numberbysection
\hybrid

\def\beq{\begin{equation}}
\def\eeq{\end{equation}}
\def\p{\partial}
\def\G{\Gamma}
\def\g{\gamma}

\def\L{{\cal L}}

\def\a{\alpha}

\def\L{{\cal L}}

\def\M{{\cal M}}

\def\H{{\cal H}}

\def\dim{{\rm dim}}

\def\Y{{\cal Y}}

\def \matrix #1 {\left(\begin{array}{cc} #1 \end{array}\right)}
\newtheorem{theo}{Theorem}[section]
\newtheorem{prop}{Proposition}[section]

\begin{document}
\begin{titlepage}
\title{A note on critical points of integrals of soliton equations}
\author{I.Krichever, \thanks{Columbia University, New York, USA and
Landau Institute for Theoretical Physics, Moscow, Russia; e-mail:
krichev@math.columbia.edu. Research is supported in part by National Science
Foundation under the grant DMS-04-05519.}}
\author{I.M.Krichever\thanks{Columbia University, New York, USA and Landau
Institute for Theoretical Physics, Moscow, Russia; e-mail:
krichev@math.columbia.edu.}
\and D.V.Zakharov \thanks{Columbia University, New York, USA; e-mail:
zakharov@math.columbia.edu}}

\date{}

\maketitle

\begin{abstract} 
We consider the problem of extending the integrals of motion of soliton equations to the space of all finite-gap solutions. We consider the critical points of these integrals on the moduli space of Riemann surfaces with marked points and jets of local coordinates. We show that the solutions of the corresponding variational problem have an explicit description in terms of real-normalized differentials on the spectral curve. Such conditions have previously appeared in a number of problems of mathematical physics.
\end{abstract}

\end{titlepage}

\section{Introduction.}
In the modern theory of integrable systems, the constants of motion of a soliton equation are defined as integrals of certain differential polynomials involving the solution. These integrals are then regarded as Hamiltonians of compatible commuting flows, and the critical points of these Hamiltonians are stationary points of the corresponding flows. The purpose of this paper is to revisit this fundamental concept of the soliton theory in a more general setting and to present unexpected connections of the critical points of the soliton integrals with physical problems, including the theory of $2D$ quantum gravity, the Whitham approach to string equations, the theory of orthogonal polynomials, and the study of Hele-Shaw flows.

First of all, it is necessary to emphasize that the Hamiltonian theory of periodic non-singular solutions of soliton equations is well-defined only for solutions of a fixed period. It is not {\it a priori} clear that the Hamiltonians considered as functionals on the space of all periodic functions with variable periods have critical points. In fact, we will show
that in the na\"ive formulation this problem has no solutions in the simplest case of the stationary KdV equation.

The problem becomes non-trivial and a posteriori interesting when extended onto the space of all finite-gap (algebro-geometric) solutions. Although these solutions are in general singular and complex functions of their arguments, their integrals can be defined in a pure algebraic-geometrical form and can be seen as functions on the moduli spaces of algebraic curves with fixed jets of local coordinates in the neighborhoods of punctures.

It is instructive to present the formulation of the extended problem and its solution in the case of the KdV equation
\beq\label{kdv} 4u_t=6uu_x-u_{xxx}.
\eeq
This equation has an infinite number of integrals of motion, which in the periodic case have the form
\beq\label{kdv-int} H_n={1\over T} \int_0^T P_n(u,u',\ldots, u^{(n+1)})\,
dx,\quad n=-1,1,3,5,\ldots,
\eeq
where the $P_n$ are certain differential polynomials in $u$:
\beq\label{kdv-int1} P_{-1}=-{u\over 2},\ \ P_1=-{u^2\over 8},\ \
P_3=-{u_x^2+2u^3\over 32}, \ldots
\eeq
Consider an arbitrary linear combination $H=H_{2n+1}+c_{2n-1}H_{2n-1}+\ldots+c_{-1}H_{-1}$ of the integrals on the space of functions with fixed period $T$, where the coefficients $c_i$ are real. The critical points of such a functional are solutions of an ordinary differential equation known as the {\it $n$-th stationary KdV equation} (\cite{nov}). For example, the stationary points of the functional
\begin{equation}
H=H_3-\frac{g_2H_{-1}}{8}
\label{Hintro}
\end{equation}
are stationary solutions of the usual KdV equation:
\beq\label{kdv1} 16\, \frac{\delta H}{\delta u}=u_{xx}-3u^2+g_2=0.
\eeq
Let us now consider $H$ as a functional on the space of all nonsingular periodic functions. Evaluating the variation of $H$ with respect to the period $T$ and setting it equal to zero, we obtain the additional condition
\beq\label{kdv3} \int_{0}^T u_x^2dx =0,
\eeq
which has no solutions when $u$ is real and smooth.

We now show how to extend the functional $H$ to the space of all periodic solutions of equation (\ref{kdv1}). The general solution of equation (\ref{kdv1}) is given by the formula
\beq\label{kdv2} u(x)=2\wp(x+x_0),
\eeq
where $\wp(x)=\wp(x|\omega,\omega')$ is the Weierstrass $\wp$-function corresponding to the elliptic curve $\Gamma$ defined by the equation
\beq\label{ellp} Y^2=E^3-g_2E-g_3=(E-E_0)(E-E_1)(E-E_3).
\eeq
The function $u$ given by equation (\ref{kdv2}) is real and smooth when the branch points $E_i$ of the elliptic curve are real, in which case we assume that $E_0<E_1<E_2$. For a generic elliptic curve, this solution is a quasi-periodic meromorphic function. For such functions, formula (\ref{kdv-int}) requires regularization.

For a real smooth solution $u(x)$ of (\ref{kdv1}) defined by an elliptic curve $\Gamma$ with real branch points,
the values of the integrals (\ref{kdv-int}) can be expressed in the following way (for details see \cite{dmn}).
There exists a unique meromorphic differential $dQ$ on $\Gamma$, called the {\it quasi-momentum differential},
with a second order pole at infinity with principal part
\beq\label{dQ1}
dQ\sim\frac{dz}{z^2}, \quad z=\frac{1}{\sqrt{E}},
\eeq
and no other singularities, and which in addition satisfies the normalization condition
\beq\label{dpnorm} \int_{E_1}^{E_2}dQ=0.
\eeq
Then the values of the integrals (\ref{kdv-int}) of a solution $u(x)$ corresponding to the curve $\Gamma$ are the coefficients of the corresponding Abelian integral $Q$ (which is the quasi-momentum of the auxiliary Sturm-Liouville operator $L=-\p_x^2+u$) at infinity:
\beq\label{Q}
Q=z^{-1}+\sum_{n=0}^\infty H_{2n-1}\, z^{2n+1}.
\eeq
For an elliptic curve with real branch points, the coefficients $H_n$ are real-valued.

As one of the authors has pointed out in \cite{non-statShrod}, the quasi-momentum differential for a generic spectral curve is uniquely defined by equation (\ref{dQ1}) and by the condition that all its periods are {\it real}:
\beq \label{norm-real}
\mbox{Im} \oint_{\gamma}dQ=0\quad\forall\gamma\in H_1(\Gamma,\mathbb{Z}).
\eeq
For elliptic curves with real branch points the two normalization conditions (\ref{dpnorm}) and (\ref{norm-real}) coincide. We can then define the functions $H_n$ on the moduli space of smooth elliptic curves as the coefficients of $dQ$ at infinity, i.e. using formula (\ref{Q}). These functions are in general complex-valued. We then define the real-valued extension of the KdV integrals (\ref{kdv-int}) to be real parts of these functions, i.e. we set
\beq\label{hdef}
\H_n=\mbox{Re}\, H_n.
\eeq

The main results of this paper are Theorems \ref{main2} and \ref{main3}, which describe the critical points of such functionals. In the case under consideration, Th. \ref{main2} implies that an elliptic curve $\Gamma$ given by equation (\ref{ellp}) is a critical point of the extended Hamiltonian $\H=\H_3$ if and only if the following condition is satisfied:
\beq\label{equa}
\mbox{Im} \oint_{\g} YdE =0 \quad\forall\gamma\in H_1(\Gamma,\mathbb{Z}).
\eeq
This equation has already appeared in the framework of the Whitham approach to string equations \cite{boutroux,kr-pen,gamba,grinevich-novikov,kapaev}.
Later similar conditions were introduced in the theory of $2D$-gravity \cite{david1,david2} and in the study of distribution of zeros in orthogonal polynomials \cite{bertola-mo}. Recently, the physical meaning of this condition in the theory of Hele-Shaw flow was revealed in \cite{wiegmann}.

Equation (\ref{equa}) was solved numerically in \cite{gamba}. It turns out that there is a unique solution in the space of smooth elliptic curves. For real $g_2>0$  the value of $g_3$ of the corresponding curve is
\beq\label{g3}
g_3={4h^2+1\over (4h^2-3)^{3/2}}\ g_2^{3/2},
\eeq
where $h\approx 3.2463822253744278875676$.

The original goal of the authors was to study the geometry of the moduli spaces of smooth curves using the extended soliton integrals as Morse functions. The connection of the variational problem for these integrals with important physical problems is to some extent an unexpected byproduct of the authors' original goal.

In the most general setting, the extended soliton integrals are defined as functions on the moduli space $\M_{g,N}({\bf n}, {\bf m})=\{(\G,P_{\alpha},[z_{\alpha}]_{\alpha},E,Q)\}$ of smooth genus $g$ algebraic curves $\G$ with $N$ marked points $P_1,\ldots,P_{N}$ with jets of local coordinates $[z_{\alpha}]$ at $P_{\alpha}$ and with a pair of Abelian integrals, denoted $E$ and $Q$ for historical reasons, having poles of orders $n_\a$ and $m_\a$ at the marked points $P_\a$. This moduli space is central to several theories with very distinct goals and origins. These include the non-linear WKB (or Whitham) theory, two-dimensional topological models, Seiberg--Witten exact solutions of $N=2$ supersymmetric gauge theories (see \cite{kr-tau,kp1,kp2} and references therein), and the theory of Hele--Shaw flows \cite{wiegmann}.  Recently, in \cite{kr-grush1,kr-grush2} it was shown that the constructions of the Whitham theory on the moduli space $\M_{g,N}^{\bf m}({\bf n})$ provide a geometrical explanation of cohomological vanishing results on the moduli spaces $\M_{g,n}$ of smooth curves with marked points.

\section{The one-point case}

The object of our study in this paper is the moduli space of Riemann surfaces with marked points and with jets of local coordinates at the marked points. Given a Riemann surface with jets of local coordinates at marked points, we construct a pair of real-normalized Abelian integrals on the surface. We then define a collection of functions on the moduli space by considering the periods of these integrals and the coefficients of their expansions at the marked points, and we consider the critical points of these functions on the leaves of a hierarchy of foliations defined on the moduli space. We show that the critical points of these functions admit an explicit analytic description. To clarify the exposition, in this section we describe the case of one marked point.

{\noindent \bf Definition.} Let $\Gamma$ be a Riemann surface, let $P$ be point on $\Gamma$. An {\it $n$-jet of local coordinates} at $P$ is an equivalence class of local coordinates at $P$, with $z$ equivalent to $z'$ if $z'=z+O(z^{n+1})$.

{\noindent \bf Definition.} Let $\omega$ be a meromorphic differential on a Riemann surface $\Gamma$. We say that $\omega$ is {\it real-normalized} if the integral of $\omega$ along any cycle in $\Gamma$ is a real number, i.e. if
\begin{equation}
\mbox{Im}\displaystyle\oint_{\gamma}\omega=0\quad\forall \gamma\in H_1(\Gamma,\mathbb{Z}).
\end{equation}

The period matrix of a Riemann surface has positive definite imaginary part, therefore a holomorphic real-normalized differential is identically zero. More generally, given a collection of principal parts with imaginary residues, there exists a unique real-normalized meromorphic differential with these principal parts:

\begin{prop} Let $\Gamma$ be a Riemann surface, let $P_1,\ldots,P_N$ be points on $\Gamma$, let $z_{\alpha}$ be local coordinates in neighborhoods of $P_{\alpha}$, and let $c_{\alpha,k}$, for $\alpha=1,\ldots, N$, $k=1,\ldots,n_{\alpha}$ be complex numbers, such that $c_{\alpha,1}$ are purely imaginary and $\sum c_{\alpha,1}=0$. Then there exists a unique real-normalized meromorphic differential $\omega$ on $\Gamma$ with poles at $P_{\alpha}$ of the form
\begin{equation}
\omega=\left(\frac{c_{\alpha,n_{\alpha}}}{z_{\alpha}^{n_{\alpha}}}+\cdots+\frac{c_{\alpha,1}}{z_{\alpha}}+
O(1)\right)dz_{\alpha}\mbox{ near }P_{\alpha},
\end{equation}
and with no other singularities.
\label{existence}
\end{prop}

{\noindent \bf Definition.} Let $\Gamma$ be a Riemann surface and let $P_0$ be a point on $\Gamma$. An {\it Abelian integral with fixed branch} $P_0$ is a multi-valued function $A(P,\gamma)$ on $\Gamma$ depending on a choice of path of integration $\gamma$ from $P_0$ to $P$, and such that its differential $dA$ is a single-valued meromorphic differential on $\Gamma$.  We say that an Abelian integral is {\it real-normalized} if the corresponding differential is.

In the neighborhood of any point that is not a residue of $dA$, the Abelian integral $A$ is locally a single-valued function. In particular, if $dA$ has no residues at $P_0$, then in a sufficiently small neighborhood of $P_0$ the Abelian integral $A$ has a selected branch corresponding to the trivial path from $P_0$.

We consider the following moduli spaces:

\begin{itemize}

\item $\M_{g,1}(n)$ is the moduli space of Riemann surfaces $\Gamma$ with a marked point $P$ and an $(n+1)$-jet of local coordinates $[z]$ at $P$. This space is our principal object of interest.

\item $\M_{g,1}^m(n)$ is the moduli space of Riemann surfaces $\Gamma$ with a marked point $P$, an $(n+1)$-jet of local coordinates $[z]$ at $P$, and a real-normalized Abelian integral $Q$ with a fixed branch at $P$, having a pole of order $m$ at $P$ and no other singularities. This space is the natural domain of the extended Hamiltonian functions that we are studying.

\item ${\bf M}_{g,1}(m,n)$ is the moduli space of Riemann surfaces $\Gamma$ with a marked point $P$, an $(n+1)$-jet of local coordinates $[z]$ at $P$, and a pair of Abelian integrals $E$ and $Q$ (not necessarily real-normalized) with fixed branches at $P$, having poles of order $n$ and $m$ at $P$, respectively, and no other singularities, where we in addition assume that the selected branch of $E$ at $P$ has the form $E=z^{-n}+O(z)$. This space is a slight modification of the moduli space introduced in the paper \cite{kp1}.

\end{itemize}

These spaces are related by the following maps:
\begin{equation}
\begin{array}{ccc}
\M_{g,1}^m(n) & \hookrightarrow & {\bf M}_{g,1}(n,m) \\
\downarrow & & \\
\M_{g,1}(n) & & \end{array}
\label{diag2}
\end{equation}
where the vertical map consists in forgetting the integrals.

Let $\G$ be a Riemann surface, and let $[z]$ be an $(n+1)$-jet of local parameters at a point $P$ on $\G$. By Prop. \ref{existence} there exists a unique real-normalized Abelian integral $E$ on $\G$ with a fixed branch at $P$, where it has a pole of order $n$ of the form $E=z^{-n}+O(z)$, and with no other singularities. There exists a unique local coordinate $z$ in the jet $[z]$ such that the fixed branch of $E$ can be written as
\begin{equation}
E=z^{-n} \mbox{ near }P,
\label{defe2}
\end{equation}
where the equality is now assumed to be exact. Therefore, the notion of real-normalized Abelian integrals gives us a natural choice of a representative in the jet $[z]$.

Let $Q$ be another Abelian integral on $\Gamma$. We can expand it in terms of this chosen local coordinate:
\begin{equation}
Q=\displaystyle\sum_{j=0}^{m} q_jz^{-j}+O(z) \mbox{ near }P.
\label{defq2}
\end{equation}
By virtue of Prop. \ref{existence} we can identify the Abelian integral $Q$ with the polynomial $q(t)=q_mt^m+\cdots+q_0$. Therefore, the vertical map in the diagram (\ref{diag2}) is a trivial bundle with fiber $\mathbb{C}^{m+1}$. Any polynomial $q(t)$ of degree $m$ determines a section of this bundle, i.e. an embedding of $\M_{g,1}(n)$ into $\M_{g,1}^m(n)$ and hence into ${\bf M}_{g,1}(n,m)$.

We now consider a collection of functions on the space ${\bf M}_{g,1}(n,m)$. Let $(\Gamma,P,[z],E,Q)$ be a point of ${\bf M}_{g,1}(n,m)$. Consider the coefficients of the expansion of the differential $QdE$ at the point $P$ with respect to the chosen local coordinate $z$:
\beq\label{qde2}
QdE=\left(\displaystyle\sum_{k=0}^{n+m}T_{k}z^{-k-1}+\displaystyle\sum_{j=1}^{\infty}H_{j}z^{j-1}\right)dz.
\eeq
The coefficients of the principal part $T_k$ act as coordinates on the space ${\bf M}_{g,1}(n,m)$. The coefficients of the regular part $H_j$ are the Hamiltonian functions that are of interest to us.

The following additional functions are only locally well defined on ${\bf M}_{g,1}(n,m)$. Fix a basis of cycles $A_i,B_i$ in $H_1(\Gamma,\mathbb{Z})$ with canonical intersection form, such that no cycle passes through $P$ or the zero-divisors of $dE$ or $dQ$. Let $\bar{\Gamma}$ denote the Riemann surface cut along these cycles and let $\partial \Gamma$ denote the boundary of
$\bar{\Gamma}$. The branches $E$ and $Q$ fixed at $P$ then extend to single-valued functions on $\bar{\Gamma}$. Consider the real-valued periods of the differentials $dE$ and $dQ$ along the chosen cycles:
\beq\label{per2}
\tau_{A_i}^{E}=\oint_{A_i}d E,\quad\tau_{B_i}^{E}=\oint_{B_i}d E,\quad
\tau_{A_i}^{Q}=\oint_{A_i}d Q,\quad\tau_{B_i}^{Q}=\oint_{B_i}d Q,\quad i=1,\ldots, g,
\eeq
and the $A$-periods of the differential $QdE$:
\beq\label{perqde2}
a_i=\oint_{A_i} QdE,\quad i=1,\ldots,g.
\eeq
The residue of the differential $QdE$ at $P$ is equal to its integral over the boundary $\p \Gamma$, which gives the only relation between these functions:
\beq\label{rel2}
2\pi iT_{0}=\sum_{i=1}^g \tau_{A_i}^Q \tau_{B_i}^E-\tau_{B_i}^Q \tau_{A_i}^E.
\eeq
Therefore, we can exclude $T_{0}$ from consideration.

The functions $T_k$ and $H_j$ do not depend on the choice of basis for $H_1(\Gamma,\mathbb{Z})$ and are well-defined on all of ${\bf M}_{g,1}(n,m)$, while the other functions defined above depend on the choice of cycles. However, in the neighborhood of any point of ${\bf M}_{g,1}(n,m)$ this choice can be made in a consistent way for all curves, therefore these functions are {\it locally} well defined. In \cite{kp1} it was shown that these functions can be used as local coordinates on ${\bf M}_{g,1}(n,m)$:

\begin{theo} Let $(\Gamma,P,[z],E,Q)$ be a point ${\bf M}_{g,1}(n,m)$ such that the zero-divisors of $dE$ and $dQ$ on $\Gamma$ do not intersect. Then in a neighborhood of this point the functions $T_{k}$ (excluding $T_0$), $\tau_{A_i}^E$, $\tau_{B_i}^E$, $\tau_{A_i}^Q$, $\tau_{B_i}^Q$ and $a_i$ (defined with respect to a consistent choice of cycles) have linearly independent differentials and therefore define a local holomorphic coordinate system for ${\bf M}_{g,1}(n,m)$.
\label{coords2}
\end{theo}

We now proceed as follows. Fix a polynomial $q(z)$ of degree $m$, and consider the embedding $\mathcal{M}_{g,1}(n)\hookrightarrow\mathcal{M}_{g,1}^m(n)\hookrightarrow{\bf M}_{g,1}(n,m)$ defined by that polynomial. We then define a hierarchy of foliations on ${\bf M}_{g,1}(n,m)$ and consider their restrictions to $\mathcal{M}_{g,1}(n)$. Finally, we define a collection of functions on $\mathcal{M}_{g,1}(n)$ and study the critical points of their restrictions to the leaves of the various foliations.

The period functions (\ref{per2}) are defined with respect to a choice of basis of $H_1(\Gamma,\mathbb{Z})$, and any two such choices are related by an element of $Sp(2n,\mathbb{Z})$. Therefore, at any point $p\in{\bf M}_{g,1}(n)$ the intersection of the null spaces of the differentials $d\tau_{A_i}^E$ and $d\tau_{B_i}^E$ in the tangent space $T_p({\bf M}_{g,1}(n))$ does not depend on the choice of basis and is therefore a well-defined subspace $T^E_p({\bf M}_{g,1}(n))$ of $T_p({\bf M}_{g,1}(n))$. Moreover, these subspaces are the tangent spaces to the level sets of some locally well-defined functions, therefore they form an integrable distribution on ${\bf M}_{g,1}(n,m)$. Similarly, the intersection of the null spaces of $d\tau_{A_i}^Q$ and $d\tau_{B_i}^Q$ is a well-defined integrable distribution $T^Q_p({\bf M}_{g,1}(n))$ on ${\bf M}_{g,1}(n)$.

The embedding $\M_{g,1}^m(n)\hookrightarrow {\bf M}_{g,1}(n,m)$ is locally defined by the condition that the functions (\ref{per2}) are real-valued. This condition does not depend on the choice of basis for $H_1(\Gamma,\mathbb{Z})$. Therefore, the distributions defined above descend to $\M_{g,1}^m(n)$. Furthermore, the space $\M_{g,1}(n)$ is a level set of the functions $T_n,\ldots,T_{n+m}$, and therefore by Prop. \ref{coords2} these distributions descend to $\M_{g,1}(n)$.

Let $\Lambda\subset\M_{g,1}(n)$ denote a joint level set of the functions $T_1,\ldots,T_{n-1}$. Let $T_p^E$ denote the intersection of the null spaces of the differentials $d\tau_{A_i}^E$ and $d\tau_{B_i}^E$ on the space $T_p(\Lambda)$. Let $\L_E\subset\Lambda$ denote an integral submanifold of this distribution. Along $\L_E$, the periods of $dE$ are constant:
\beq\label{leafE2}
\delta \oint_\g dE=0\quad\forall\delta\in T_p(\L_E)=T_p^E\mbox{ and }\forall\gamma\in H_1(\Gamma,\mathbb{Z}).
\eeq
Similarly, let $T_p^Q\subset T_p(\Lambda)$ denote the directions along which the periods of $dQ$ do not change, and let  $\L_Q\subset\Lambda$ denote the corresponding integral submanifold:
\beq\label{leafQ2}
\delta \oint_\g dQ=0\quad\forall\delta\in T_p(\L_Q)=T_p^Q\mbox{ and }\forall\gamma\in H_1(\Gamma,\mathbb{Z}).
\eeq
Finally, let $\L=\L_E\cap\L_Q$ denote a submanifold along which the periods of both $dE$ and $dQ$ are constant:
\beq\label{leaf2}
\delta \oint_\g dE=\delta\oint_\g dQ=0\quad\forall\delta\in T_p(\L)=T_p^E\cap T_p^Q\mbox{ and }\forall\gamma\in H_1(\Gamma,\mathbb{Z}).
\eeq
We now define a real-valued function on $\mathcal{M}_{g,1}(n)$ and study its critical points on the leaves of the foliations defined above. Let $c_j$, $d_j$, for $j=1,\ldots,M$ be real coefficients. Consider the following function on $\mathcal{M}_{g,1}(n)$:
\beq\label{ham2}
\H=\sum_{j=1}^{M}\left( c_{j}\mbox{Re} H_{j}+d_j\mbox{Im}H_j\right).
\eeq
The following theorem, which is the principal result of this chapter, describes the critical points of $\H$ restricted to the leaves of the foliations defined above:

\begin{theo} Let $(\Gamma,P,[z])$ be a point in $\mathcal{M}_{g,1}(n)$, which we assume to be embedded in $\mathcal{M}_{g,1}^m(n)$ via the polynomial $q(z)$. Let $E$ and $Q$ be the real-normalized Abelian integrals defined by equations (\ref{defe2}) and (\ref{defq2}), and let $D$ denote the intersection of the zero-divisors of $dE$ and $dQ$.

\noindent {\bf A.} The point $(\Gamma,P,[z])$ is a critical point of $\H$ restricted to $\L$ if and only if

(i) There exists a meromorphic function $Y$ on $\Gamma$ with a pole at $P$ with the principal part
\beq\label{th221}
Y=\sum_{j=1}^{\infty} (c_{j}-id_j)z^{-j}+O(z)\mbox{ near }P,
\eeq
with poles at the divisor $D$, and no other singularities.

\noindent {\bf B.} The point $(\Gamma,P,[z])$ is a critical point of $\H$ restricted to $\L_E$ if and only if there exists a function $Y$ satisfying condition (i) and in addition the following condition:

(ii) the differential $YdE$ is real-normalized, i.e.
\beq\label{th222}
\mbox{Im} \oint_{\g} YdE=0\quad\forall\gamma\in H_1(\Gamma,\mathbb{Z}).
\eeq

\noindent {\bf C.} The point $(\Gamma,P,[z])$ is a critical point of $\H$ restricted to $\L_Q$ if and only if there exists a function $Y$ satisfying condition (i) and in addition the following condition:

(iii) the differential $YdQ$ is real-normalized, i.e.
\beq\label{th223}
\mbox{Im}\oint_{\g} YdQ=0\quad\forall\gamma\in H_1(\Gamma,\mathbb{Z}).
\eeq

\noindent {\bf D.} The point $(\Gamma,P,[z])$ is a critical point of $\H$ restricted to $\Lambda$ if and only if there exists a function satisfying conditions (i), (ii) and (iii).
\label{main2}
\end{theo}

{\bf Remark.} The correspondence between critical points of integrals of soliton equations and meromorphic functions with prescribed singularities on the spectral curve, in other words part A of our theorem, has long been a central notion in the theory of finite-gap integration (see, for example, the survey \cite{dmn}). However, this result has never before been stated in such an explicit form. Moreover, the known results only apply to the case when the equations are Hamiltonian, which corresponds to the case when the zero-divisors of the differentials $dE$ and $dQ$ do not intersect.

The Hamiltonian structure of the integrable systems associated with the moduli space $\M_{g,N}(n)$ is determined by the Novikov--Veselov--Witten--Seiberg symplectic form, which is defined as follows. Consider a leaf $\L$ of the foliation defined above, anc consider the fibration $\mathcal{N}^g$ over $\L$ whose fiber over a point $(\Gamma,P,[z])$ is the $g$-th symmetric power of $\Gamma$. This fibration has dimension $2g$ and is the phase space for the corresponding integrable system. The Novikov--Veselov--Witten--Seiberg symplectic form is defined as
\begin{equation}
\omega_{\M}=\displaystyle\sum_{i=1}^g da_i \wedge d\omega_i,
\end{equation}
where $\omega_i$ are normalized holomorphic differentials on $\Gamma$.

If $(\Gamma,P,[z])$ is a point at which the zero-divisors of $dE$ and $dQ$ intersect, then according to Th. \ref{coords2} the coordinates (\ref{perqde2}) are not independent on the leaf $\L$. Therefore, the symplectic form $\omega_{\M}$ is degenerate at $(\Gamma,P,[z])$, and therefore the soliton equation is not Hamiltonian.

\bigskip

{\bf Example: elliptic solutions of KdV.} Before proceeding with the proof of the theorem, we show how it applies to the example considered in the introduction. Assume that $n=2$ and $m=1$, and let $g$ for now be arbitrary. Consider the embedding $\mathcal{M}_{g,1}(2)\hookrightarrow\mathcal{M}_{g,1}^1(2)$ given by the polynomial $q(z)=z$. Suppose furthermore that we are on a leaf $\mathcal{L}_E$ along which the periods of $dE$ are identically zero, so that the Abelian integral $E$ is actually a single-valued function of degree two. Thus, at every point $(\Gamma,P,[z])\in \mathcal{L}_E$, $\Gamma$ is a hyperelliptic curve of genus $g$ and the Abelian integral $E$ is a single-valued function on $\Gamma$ of degree two, having a second order pole at $P$ of the form $E=z^{-2}$ and no other singularities.

The differential $QdE$ is odd with respect to the hyperelliptic involution and therefore has the following expansion at $P$:
\begin{equation}
QdE=-2\left(z^{-4}+T_1 z^{-2}+H_1+H_3 z^{2}+O(z^4)\right) dz.
\end{equation}
We consider the critical points of the function $\mathcal{H}=\mbox{Re } H_3$, which are the same as those of the function (\ref{Hintro}), since $H_{-1}=T_1$ is constant along $\mathcal{L}_E$. According to the theorem, a point $(\Gamma,P,[z])$ is a critical point of $H$ restricted to $\L$ if there exists a meromorphic function $Y$ on $\Gamma$ with a triple pole at $P$ of the form
\begin{equation}
Y=z^{-3}+O(z)\mbox{ near }P.
\end{equation}
This is only possible if $g=1$, so the function $\mathcal{H}$ has a critical point only on the leaf corresponding to elliptic curves. The function $Y$ then satisfies an equation of the form (\ref{ellp}). Therefore, $(\Gamma,P,[z],E,Y,Q)$ are the spectral data of a stationary genus one solution of the KdV equation. If the point $(\Gamma,P,[z])$ is a critical point of $\mathcal{H}$ on all of $\L_E$, then the function $Y$ satisfies the condition
\begin{equation}
\mbox{Im} \oint_{\g} YdE=0\quad\forall\gamma\in H_1(\Gamma,\mathbb{Z}).
\end{equation}

The main tool in the proof of the theorem is a homomorphism
\beq\label{dmap2}
f: T_p(\Lambda)\rightarrow \Omega(\Gamma,P,E,Q)
\eeq
from the tangent space to $\Lambda$ at a point $p=(\Gamma,P,[z])$ to a space of certain multi-valued differentials on $\Gamma$. This homomorphism was introduced in Whitham theory in \cite{kr-tau} (see more details in \cite{kp1,kp2}) and modified for the case real normalized differentials in \cite{dzhamay}. It is defined as follows.

Let $\bar{\Gamma}$ denote the surface cut along a choice of basis cycles, and let $\sum \nu_s E_s$ be the zero divisor of the differential $dE$. On $\bar{\Gamma}$, the branch of $E$ fixed at $P$ can be used as a local coordinate in the neighborhood of any point except $E_s$. This allows us to lift a tangent vector $\delta\in T_p(\Lambda)$ to a tangent vector on the total space of the universal curve over $\L_E$ and defines the partial derivative of $Q=Q(E)$ for fixed levels of $E$, which we denote by $\delta Q|_E$.

The derivative $\delta Q|_E$ is holomorphic away from the points $E_s$, where it has poles of order $\nu_s$, and from the cuts. Along the two sides of a cut $A_i$ or $B_i$, the values of $E$ and $Q$ have jumps equal to the corresponding periods of $dE$ and $dQ$:
\begin{equation}
Q(E+\tau_{A_i}^E)-Q(E)=\tau_{A_i}^Q,\quad Q(E+\tau_{B_i}^E)-Q(E)=\tau_{B_i}^Q.
\end{equation}
Therefore, the derivative $\delta Q|_E$ has jumps along the cuts $A_i$ and $B_i$ equal to
\begin{equation}
\delta Q(E+\tau_{A_i}^E)-\delta Q(E)=\delta \tau_{A_i}^Q-\frac{dQ}{dE}\delta\tau_{A_i}^E,\quad
\delta Q(E+\tau_{B_i}^E)-\delta Q(E)=\delta \tau_{B_i}^Q-\frac{dQ}{dE}\delta\tau_{B_i}^E.
\end{equation}

The map $f$ associates to a tangent vector $\delta$ the multi-valued meromorphic differential
\beq\label{diff}
f(\delta)=\Omega_{\delta}=(\delta Q|_E )dE.
\eeq
Since all the higher order terms in the expansion (\ref{qde2}) of $QdE$ at $P$ are fixed on $\Lambda$, the differential $\Omega_{\delta}$ has at most a simple pole at $P$. At the points $E_s$, the zeroes of $dE$ cancel the poles of $\delta Q|_E$, hence $\Omega_{\delta}$ has no other singularities. Finally, along the cuts $A_i$ and $B_i$ the differential $\Omega_{\delta}$ has jumps which are linear combinations of the differentials $dE$ and $dQ$ with real coefficients.

The space $\Omega(\Gamma,P,E,Q)$ of differentials with these properties has real dimension $6g$, the same as $\L$. It is clear from the construction that the image of a tangent vector lying inside one of the subspaces $T_p(\L_E)$ or $T_p(\L_Q)$ is a differential with jumps proportional only to $dE$ or $dQ$, respectively, and that the image of their intersection $T_p(\L)=T_p(\L_E)\cap T_p(\L_Q)$ lies inside the space of holomorphic differentials $\Omega^h(\Gamma)\subset\Omega(\Gamma,P,E,Q)$. A complete description of the homomorphism (\ref{dmap2}) was given in \cite{kp1}:
\begin{prop}
Let $D=\sum\mu_sQ_s$ denote the intersection of the zero divisors of $dE$ and $dQ$. Then the image of the tangent space $T_p(\Lambda)$ under the homomorphism $f$ of (\ref{dmap2}) is the subspace $\Omega_D(\Gamma,P,E,Q)\subset\Omega(\Gamma,P,E,Q)$ of differentials having zeroes at $Q_s$ of orders $\mu_s$. In particular, if the zero divisors of $dE$ and $dQ$ do not intersect, then the map (\ref{dmap2}) is an isomorphism.
\label{map2}
\end{prop}
The space $\Omega(\Gamma,P,E,Q)$ is spanned over $\mathbb R$ by the following collection of differentials:
\begin{itemize}

\item $2g$ single-valued holomorphic differentials $\Omega_{A_i}^h, \Omega_{B_i}^h, i=1,\ldots,g$ normalized by the conditions
\beq\label{normah}
 \mbox{Im} \left(\oint_{A_j}\Omega_{A_i}^h\right)=\delta_{ij},\quad
\mbox{Im} \left(\oint_{B_j}\Omega_{A_i}^h\right)=0,\quad j=1,\ldots,g,
\eeq
\beq\label{normbh}
\mbox{Im} \left(\oint_{A_j}\Omega_{B_i}^h\right)=0,\quad
\mbox{Im} \left(\oint_{B_j}\Omega_{B_i}^h\right)=\delta_{ij},\quad j=1,\ldots,g.
\eeq
\item $2g$ differentials $\Omega_{A_i}^E, \Omega_{B_i}^E, i=1,\ldots,g$ having a simple pole at $P$ and a single jump equal to $dE$ along the cuts $A_i$ or $B_i$, respectively:
\beq\label{jumpe}
\left(\Omega_{A_i}^E\right)^+-\left(\Omega_{A_i}^E\right)^-=dE\mbox{ along }A_i,\quad
\left(\Omega_{B_i}^E\right)^+-\left(\Omega_{B_i}^E\right)^-=dE\mbox{ along }B_i,
\eeq
normalized by the conditions
\beq\label{norme}
\mbox{Im}\oint_{A_j} \Omega_{A_i}^E=0,\quad
\mbox{Im}\oint_{B_j} \Omega_{A_i}^E=0,\quad
\mbox{Im}\oint_{A_j} \Omega_{B_i}^E=0,\quad
\mbox{Im}\oint_{B_j} \Omega_{B_i}^E=0,\quad j=1,\ldots,g.
\eeq
\item $2g$ differentials $\Omega_{A_i}^Q, \Omega_{B_i}^Q, i=1,\ldots,g$ having a simple pole at $P$ and a single jump equal to $dQ$ along the cuts $A_i$ or $B_i$, respectively:
\beq\label{jumpq}
\left(\Omega_{A_i}^Q\right)^+-\left(\Omega_{A_i}^Q\right)^-=dQ\mbox{ along }A_i,\quad
\left(\Omega_{B_i}^Q\right)^+-\left(\Omega_{B_i}^Q\right)^-=dQ\mbox{ along }B_i,
\eeq
normalized by the conditions
\beq\label{normq}
\mbox{Im}\oint_{A_j} \Omega_{A_i}^Q=0,\quad
\mbox{Im}\oint_{B_j} \Omega_{A_i}^Q=0,\quad
\mbox{Im}\oint_{A_j} \Omega_{B_i}^Q=0,\quad
\mbox{Im}\oint_{B_j} \Omega_{B_i}^Q=0,\quad j=1,\ldots,g.
\eeq
\end{itemize}

{\bf Remark.} The integral of a multi-valued differential along a curve in $\bar{\Gamma}$ is not determined by the homology class of that curve in $\Gamma$. For this reason, the integrals in the above formulas are assumed to be taken over the particular choice of $A_i$ and $B_i$ that we have made.

{\bf Proof of Theorem \ref{main2}} We first prove the statement of the theorem for the critical points of $\H$ restricted to the open subset of $\Lambda$ where the differentials $dQ$ and $dE$ have no common zeros.

Given a set of coefficients $c_j, d_j$ there is a unique real-normalized meromorphic differential $dY$ of the second kind whose Abelian integral has a single pole at $P$ of the form (\ref{th221}). By the definition of $\H$ the derivative $\delta \H$ along a tangent vector $\delta$ is equal
to
\beq\label{var2}
\delta \H=\mbox{Re}\,\mbox{Res}_{P} \left(Y\delta Q|_EdE\right)=
-{1\over 2\pi}\mbox{Im}\left(\oint_{\p \Gamma}Y\Omega_{\delta}\right),
\eeq
where $\p \Gamma$ is the boundary of $\bar{\Gamma}$. Hence, the point $(\Gamma,P,[z])$ is a critical point of $\H$ restricted to $\L$ if and only if
\beq\label{var2'}
\mbox{Im}\left(\oint_{\p \Gamma}Y\Omega\right)=\displaystyle\sum_{i=1}^g\left(\oint_{A_i}dY
 \cdot\mbox{Im}\oint_{B_i}\Omega-\oint_{B_i}dY\cdot\mbox{Im}\oint_{A_i}\Omega\right)=0
\eeq
for all holomorphic differentials $\Omega\in \Omega^h(\Gamma)$. Plugging the differentials $\Omega_{A_i}^h$ and $\Omega_{B_i}^h$ in this equation, we see that all the periods of $Y$ are zero:
\beq\label{imy2}
\oint_\g dY=0\quad\forall\gamma\in H_1(\G,\mathbb{Z}).
\eeq
Therefore $Y$ is a single-valued meromorphic function, which proves part A of the theorem. To prove parts B and C of the theorem, we plug the remaining differentials $\Omega_{A_i}^E$, $\Omega_{B_i}^E$, $\Omega_{A_i}^Q$, $\Omega_{B_i}^Q$ into equation (\ref{var2'}), which gives us conditions (\ref{th222}) and (\ref{th223}). Finally, since the tangent space $T_p(\Lambda)$ is spanned by $T_p(\L_E)$ and $T_p(\L_Q)$, part D follows from parts B and C.

The proof of the theorem in the general case goes along the same lines. Suppose that $p=(\Gamma,P,[z])$ is a critical point of $\H$ restricted to $\L$, and let
\beq\label{divis}
D=\sum_s \mu_s Q_s=(dE)_0\cap (dQ)_0.
\eeq
be the intersection of the zero-divisors of $dE$ and $dQ$. The image of the tangent subspace $T_p(\L)$ under $f$ is the space of holomorphic differentials on $\Gamma$ vanishing at $D$, which we denote $\Omega^h_D(\Gamma)$. Let $\Y_{D,c_i}$ denote space of real-normalized Abelian integrals on $\Gamma$ having poles of orders $\mu_s$ at $Q_s$ and a pole at the marked point $P$ of the form (\ref{th221}). This is an affine space of real dimension $\dim_{\mathbb{R}} \,\Y_{D,c_i}=2\sum_s \mu_s=2d$. The period vectors of the Abelian integrals in $\Y_{D,c_i}$ form an affine subspace $\Pi_{D,c_i}$ of $\mathbb{R}^{2g}$ of dimension $2d-2h^0(\Gamma,D)+2$.

It is easy to see that we can evaluate the derivative $\delta \H$ by the same formula (\ref{var2}) using any $Y\in\Y_{D ,c_i}$. Therefore, equation (\ref{var2'}) holds for all $Y\in\Y_{D,c_i}$ and for all $\Omega\in \Omega^h_D(\Gamma)$. This gives $\dim_{\mathbb{R}}\,\Omega^h_D(\Gamma)=2h^1(\Gamma,D)=2g-2d+2h^0(\Gamma,D)-2$ linearly independent conditions on the elements of $\Pi_{D,c_i}$. Therefore, $\Pi_{D,c_i}$ is actually a linear subspace of $\mathbb{R}^{2n}$. Therefore, there exists an Abelian integral $Y\in\Y_{D,c_i}$ with trivial periods, that is to say a single-valued meromorphic function, which proves part A of the theorem.

Let $Y_{D,c_i}$ denote the subspace of single-valued functions in $\Y_{D,c_i}$, we have shown above that $\dim_{\mathbb{R}}\, Y_{D,c_i}=2h^0(\Gamma,D)-2$. The image of $T_p(\L_E)$ under $f$ is the space of differentials with jumps proportional to $dE$ with real coefficients, having a simple pole at $P$ and zeroes at the divisor $D$. We denote this space by $\Omega_D(\Gamma,P,E)$. Let $J_e:\Omega_D(\Gamma,P,E) \rightarrow \mathbb{R}^{2g}$ be the map which associates to each differential its jumps along the basis cycles. Since $\dim_{\mathbb{R}}\, \Omega_D(\Gamma,P,E)=4g-2d$ and $\dim_{\mathbb{R}}\, \mbox{Ker}\, J_E=\dim_{\mathbb{R}}\, \Omega^h_D(\Gamma)=2h^1(\Gamma,D)$, we see that $\dim_{\mathbb{R}}\,\mbox{Im} \,J_E=2g+2-2h^0(\Gamma,D)=2g-\dim_{\mathbb{R}} \,Y_{D,c_i}$. Therefore, formula (\ref{var2'}) gives $2g-\dim_{\mathbb{R}}\, Y_{D,c_i}$ linear relations on the $2g$ imaginary parts of the periods of $YdE$. Therefore, there exists a $Y\in Y_{D,c_i}$ such that the differential $YdE$ is real-normalized, which proves part B. The proof of part C is identical.

Finally, to prove part D, consider the image $\Omega_D(\Gamma,P,E,Q)$ of $T_p(\Lambda)$ under $f$, and let $J_{E,Q}:\Omega_D(\Gamma,P,E,Q)\rightarrow \mathbb{R}^{2g}\oplus \mathbb{R}^{2g}$ be the map that associates to each differential its jumps, which are a real linear combination of $dE$ and $dQ$. Since $\dim_{\mathbb{R}}\,\Omega_D(\Gamma,P,E,Q)=6g-2d$ and $\dim_{\mathbb{R}}\,\mbox{Ker}\, J_{E,Q}=\dim_{\mathbb{R}}\,\Omega^h_D(\Gamma)=2h^1(\Gamma,D)$, formula (\ref{var2'}) gives $\dim_{\mathbb{R}}\,\mbox{Im} \, J_{E,Q}=4g+2-2h^0(\Gamma,D)$ linear conditions on the $4g$ imaginary parts of the periods of $YdE$ and $YdQ$. Therefore, there exists a $Y\in Y_{D,c_i}$ such that both $YdE$ and $YdQ$ are real-normalized, which completes the proof of the theorem.

\section{The $N$-point case}

In this section, we extend our problem to the moduli space of curves with several marked points and with jets of local coordinates at those points. Given a curve with jets of local coordinates, we define a real-normalized Abelian integral with poles determined by those jets. The principal difference in the $N$-point case is that is that we can only fix a branch of an Abelian integral at one marked point. At all other marked points our Abelian integrals will be truly multi-valued functions. This means that we cannot use the integral $E$ to specify a unique local coordinate in the jet, as we did in (\ref{defe2}), and therefore the functions (\ref{qde2}) are no longer globally defined on the moduli space. Furthermore, the expansion (\ref{defq2}) of the Abelian integral $Q$ depends on the choice of the local coordinate, therefore the map forgetting the second Abelian integral $Q$ is in general no longer a trivial bundle. In this section, we show that under certain assumptions the imaginary parts of the functions (\ref{qde2}) are well-defined, which allows us to construct a somewhat more limited extension of the results of the previous chapter.

Let ${\bf n}=(n_1,\ldots,n_N)$ and ${\bf m}=(m_1,\ldots,m_N)$ denote a pair of multi-indices. We assume that $n_{\alpha}>0$ and $m_{\alpha}\geq 0$ for all $\alpha$, and that $m_{\alpha}\leq n_{\alpha}$ for $\alpha=2,\ldots,N$; the meaning of the last condition will become clear later. We consider the following spaces:

\begin{itemize}

\item $\M_{g,N}({\bf n})$ is the moduli space of Riemann surfaces $\Gamma$ with marked points $P_1,\ldots,P_N$, with an $(n_1+1)$-jet of local coordinates $[z_1]$ at $P_1$ and with $n_{\alpha}$-jets of local coordinates $[z_{\alpha}]$ at $P_{\alpha}$ for $\alpha=2,\ldots,n$.

\item $\M_{g,N}^{\bf m}({\bf n})$ is the moduli space of Riemann surfaces $\Gamma$ with marked points $P_1,\ldots,P_N$, with an $(n_1+1)$-jet of local coordinates $[z_1]$ at $P_1$ and with $n_{\alpha}$-jets of local coordinates $[z_{\alpha}]$ at $P_{\alpha}$ for $\alpha=2,\ldots,N$, and with a real-normalized Abelian integral $Q$ with fixed branch at $P_1$, having poles of orders $m_{\alpha}$ at $P_{\alpha}$ and no other singularities.

\item ${\bf M}_{g,N}({\bf n},{\bf m})$ is the moduli space of Riemann surfaces $\Gamma$ with marked points $P_1,\ldots,P_N$, with an $(n_1+1)$-jet of local coordinates $[z_1]$ at $P_1$ and with $n_{\alpha}$-jets of local coordinates $[z_{\alpha}]$ at $P_{\alpha}$ for $\alpha=2,\ldots,n$, and with a pair of Abelian integrals $E$ and $Q$ (not necessarily real-normalized) with fixed branches at $P_1$, having poles of orders $n_{\alpha}$ and $m_{\alpha}$ at $P_{\alpha}$, respectively, and no other singularities, where we in addition assume that $E=z_1^{-n_1}+O(z_1)$ near $P_1$ and $E=z_{\alpha}^{-n_{\alpha}}+O(1)$ near $P_{\alpha}$ for $\alpha=2,\ldots,N$. Here we assume for simplicity that the integral $E$ has no logarithmic singularities.

\end{itemize}

These spaces are related by the following maps:
\begin{equation}
\begin{array}{ccc}
\M_{g,N}^{\bf m}({\bf n}) & \hookrightarrow & {\bf M}_{g,N}({\bf n},{\bf m})\\
\downarrow & & \\
\M_{g,N}({\bf n}) & &\end{array}
\label{diagram3}
\end{equation}
where the vertical map consists in forgetting the integrals.

Let $[z_1]$ be a $(n_1+1)$-jet of local coordinates at $P_1$ and let $[z_{\alpha}]$ be $n_{\alpha}$-jets of local coordinates at $P_{\alpha}$ for $\alpha=2,\ldots,N$. According to Prop. \ref{existence} there exists a unique real-normalized Abelian integral $E$ on $\Gamma$ with a fixed branch at $P_1$ and with poles of orders $n_{\alpha}$ at $P_{\alpha}$ of the form
\begin{equation}
E=z_{1}^{-n_1}+O(z_1)\mbox{ near }P_1,\quad E=z_{\alpha}^{-n_{\alpha}}+O(1)\mbox{ near }P_{\alpha},\quad \alpha=2,\ldots,N,
\label{defe3}
\end{equation}
and with no other singularities. There exists a unique local coordinate $z_1$ in the jet $[z_1]$ such that $E=z_1^{-n_1}$ exactly. However, at the other points $P_{\alpha}$ the Abelian integral $E$ does not have a selected branch, so we cannot choose a representative of the jet $[z_{\alpha}]$ in this way.

Let $E'$ and $E''$ denote two branches of the Abelian integral $E$ at a point $P_{\alpha}$. We have $E''=E'+C$, where the constant $C$ is real if the integral $E$ is real-normalized. Let $z'_{\alpha}$ and $z''_{\alpha}$ denote the unique representatives of the jet $[z_{\alpha}]$ such that $E'=(z'_{\alpha})^{-n_{\alpha}}$ and $E''=(z''_{\alpha})^{-n_{\alpha}}$ near $P_{\alpha}$. It is easy to check that $z'_{\alpha}$ and $z''_{\alpha}$ are related as follows:
\begin{equation}
z''_{\alpha}=z'_{\alpha}+\displaystyle\sum_{j=1}^{\infty}C_j (z'_{\alpha})^{jn_{\alpha}+1},
\label{ambig}
\end{equation}
where there the constants $C_j$ are polynomials in $C$ with rational coefficients.

Now let $Q$ be another Abelian integral on $\Gamma$. At the point $P_1$ it can be expanded in terms of the chosen local coordinate $z_1$:
\begin{equation}
Q=\displaystyle\sum_{j=0}^{m_1} q_{1,j}z_1^{-j}+R_{1}^Q\log z_1+O(z_1)\mbox{ near }P_1.
\label{defq31}
\end{equation}
To expand $Q$ at the other points $P_{\alpha}$, we need to choose local coordinates at those points. Let $z_{\alpha}$ denote a choice of local coordinate in the jet $[z_{\alpha}]$ corresponding to some branch of $E$. Then $Q$ can be expanded at $P_{\alpha}$:
\begin{equation}
Q=\displaystyle\sum_{j=1}^{m_{\alpha}} q_{\alpha,j}z_{\alpha}^{-j}+R_{\alpha}^Q \log z_{\alpha}+O(1)
\mbox{ near }P_{\alpha},\quad \alpha=2,\ldots,N.
\label{defq3a}
\end{equation}
If we now assume that $m_{\alpha}\leq n_{\alpha}$, then changing the representative of $[z_{\alpha}]$ according to (\ref{ambig}) does not affect the principal terms in the expansion of $Q$. Therefore, if $m_{\alpha}\leq n_{\alpha}$ for $\alpha=2,\ldots,N$, then the map $\M_{g,N}^m(n)\rightarrow\M_{g,N}(n)$ is a trivial fiber bundle with fiber $\mathbb{C}^{m_1+1}\oplus \mathbb{C}^{m_2}\oplus\cdots\oplus\mathbb{C}^{m_N}\oplus \mathbb{R}^{N-1}$. To define an embedding $\M_{g,N}({\bf n})$ into $\M_{g,N}^{\bf m}({\bf n})$, we need to choose polynomials $q_{\alpha}(t)$ of degrees $m_{\alpha}$, such that $q_{\alpha}(t)$ have no constant terms for $\alpha=2,\ldots,N$, and $N$ imaginary numbers $R_{\alpha}^Q$ adding up to zero. If we do not make the assumption that $m_{\alpha}\leq n_{\alpha}$ for $\alpha=2,\ldots,N$, then the map is only locally trivial, and we cannot define a global embedding $\M_{g,N}({\bf n})\hookrightarrow\M_{g,N}^{\bf m}({\bf n})$.

We now define a collection of functions on ${\bf M}_{g,N}({\bf n},{\bf m})$. First, we consider the residues of the differential $dQ$:
\begin{equation}
\label{res3}
R^Q_{\alpha}=\mbox{Res}_{P_{\alpha}}dQ, \quad \alpha=1,\ldots,N.
\end{equation}
These functions satisfy the relation $\sum R^Q_{\alpha}=0$.

The remaining functions are only locally well-defined on ${\bf M}_{g,N}({\bf n},{\bf m})$. Fix a basis of cycles $A_i,B_i$ in $H_1(\Gamma,\mathbb{Z})$ with canonical intersection form, such that no cycle passes through any of the $P_i$ or the zero-divisors of $dE$ and $dQ$. Let $\bar{\Gamma}$ denote the Riemann surface cut along these cycles and let $\partial \Gamma$ denote the boundary of $\bar{\Gamma}$. Fix also a choice of paths $\gamma_{2},\ldots,\gamma_{n}$ from $P_1$ to $P_2,\ldots,P_n$ not intersecting any of the basis cycles. The fixed branches of $E$ and $Q$ then extend to functions on $\bar{\Gamma}\,\backslash\!\displaystyle\cup \gamma_i$.

Let $z_{\alpha}$ be the unique local coordinate in the jet $[z_{\alpha}]$ such that for the chosen branch of $E$ we have $E=z_{\alpha}^{-n_{\alpha}}$ near $P_{\alpha}$. We consider the coefficients of the expansion of $QdE$ at the marked points, with respect to the chosen local coordinates $z_{\alpha}$:
\begin{equation}
QdE=\left(\displaystyle\sum_{k=0}^{n_{\alpha}+m_{\alpha}}T_{\alpha,k}z_{\alpha}^{-k-1}-n_{\alpha}R_{\alpha}^Qz^{-n_{\alpha}-1} \log z_{\alpha}+\displaystyle\sum_{j=1}^{\infty}H_{\alpha,j} z_{\alpha}^{j-1}\right)d z_{\alpha}.
\label{qde3}
\end{equation}
We also consider the periods of the differentials $dE$ and $dQ$ along the chosen cycles:
\beq\label{per3}
\tau_{A_i}^{E}=\oint_{A_i}d E,\quad\tau_{B_i}^{E}=\oint_{B_i}d E,\quad
\tau_{A_i}^{Q}=\oint_{A_i}d Q,\quad\tau_{B_i}^{Q}=\oint_{B_i}d Q,\quad i=1,\ldots, g,
\eeq
and the $A$-periods of the differential $QdE$:
\beq\label{perqde3}
a_i=\oint_{A_i} QdE,\quad i=1,\ldots,g.
\eeq
The sum of the residues of the differential $QdE$ is equal to its integral over $\partial\Gamma$. That gives the only relation between these functions:
\beq\label{rel3}
2\pi i\sum_{\a=1}^N T_{\a,0}=\sum_{i=1}^g\left( \tau_{A,i}^Q \tau_{B,i}^E-\tau_{B,i}^Q \tau_{A,i}^E\right),
\eeq
so therefore we exclude $T_{1,0}$ from consideration.

As before, these functions are locally well defined on ${\bf M}_{g,N}({\bf n},{\bf m})$ and have well-defined differentials, and they can be used as holomorphic coordinates on ${\bf M}_{g,N}({\bf n},{\bf m})$. The corresponding generalization of Theorem \ref{coords2} was also proved in \cite{kp1}:

\begin{theo} Let $P=(\Gamma,P_{\alpha},[z_{\alpha}],E,Q)$ be a point on ${\bf M}_{g,N}({\bf n},{\bf m})$ such that the zero-divisors of $dE$ and $dQ$ on $\Gamma$ do not intersect. Then in a neighborhood of $P$ the collection of $5g+2N-2+\sum(n_{\alpha}+m_{\alpha})$ functions $R^Q_{\alpha}$, $T_{\alpha,k}$ (excluding $R^Q_1$ and $T_{1,0}$), $\tau_{A_i}^E$, $\tau_{B_i}^E$, $\tau_{A_i}^Q$, $\tau_{B_i}^Q$ and $a_i$ have linearly independent differentials and thus define a local holomorphic coordinate system for ${\bf M}_{g,N}({\bf n},{\bf m})$.
\label{coords3}
\end{theo}

In the $N$-point case, the local coordinates $z_{\alpha}$ and hence the functions $H_{\alpha,k}$ are no longer globally well-defined on ${\bf M}_{g,N}({\bf n},{\bf m})$ for $\alpha=2,\ldots,N$. Hence, we cannot talk of their critical points.
However, on the smaller subspace $\M_{g,N}^{\bf m}({\bf n})$ the imaginary parts of these functions are well-defined. Indeed, any two choices of local coordinate $z_{\alpha}$ in the jet $[z_{\alpha}]$ are related by equation (\ref{ambig}). On $\M_{g,N}^{\bf m}({\bf n})$ the coefficients $C_j$ are real, and therefore substituting them into expansion (\ref{qde3}) does not change the imaginary parts of the functions $H_{\alpha,k}$.

It is now clear how to generalize the results of the previous section to the $N$-point case. Fix $N$ polynomials $q_{\alpha}(t)$ of degrees $m_{\alpha}$ with no constant terms for $\alpha=2,\ldots,N$, and $N$ imaginary numbers $R^Q_{\alpha}$ adding up to zero. Consider the embedding $\M_{g,N}({\bf n})\hookrightarrow\M_{g,N}^{\bf m}({\bf n})$ defined by these data.

We introduce a hierarchy of foliations on $\M_{g,N}({\bf n})$.  Let $\Lambda'\subset\M_{g,N}({\bf n})$ be a submanifold along which the functions $T_{\alpha,k}$ for $k=1,\ldots,m_{\alpha}+n_{\alpha}$ are constant:
\beq\label{leafL'3}
\delta T_{\alpha,k}=0\quad\forall k>0\mbox{ and }\forall\delta\in T(\Lambda').
\eeq
Let $\Lambda\subset\Lambda'$ denote a submanifold along which the functions $T_{\alpha,k}$ for all $k=0,\ldots,m_{\alpha}+n_{\alpha}$ are constant:
\beq\label{leafL3}
\delta T_{\alpha,k}=0\quad\forall k\geq 0\mbox{ and }\forall \delta\in T(\Lambda).
\eeq
Let $\L_E\subset\Lambda$ denote a submanifold along which the periods of $dE$ are constant:
\beq\label{leafE3}
\delta \oint_\g dE=0\quad\forall\delta\in T(\L_E)\mbox{ and }\forall\gamma\in H_1(\Gamma,\mathbb{Z}),
\eeq
let $\L_Q \subset\Lambda$ denote a submanifold along which the periods of $dQ$ are constant:
\beq\label{leafQ3}
\delta\oint_\g dQ=0\quad\forall\delta\in T(\L_Q)\mbox{ and }\forall\gamma\in H_1(\Gamma,\mathbb{Z}).
\eeq
Finally, let $\L=\L_E\cap\L_Q$ denote a submanifold along which the periods of $dE$ and $dQ$ are both constant:
\beq\label{leaf3}
\delta \oint_\g dE=\delta\oint_\g dQ=0\quad\forall\delta\in T(\L)\mbox{ and }\forall\gamma\in H_1(\Gamma,\mathbb{Z}).
\eeq
We now define a real-valued function on $\M_{g,N}({\bf n})$ and study its critical points on the leaves of the foliations described above. Let $d_{\a,j}$ for $\alpha=1,\ldots,N$, $j=1,\ldots,M_{\alpha}$ be a finite collection of real constants. As we have noted above, the function
\beq\label{ham3}
\H=\sum_{\a,j}d_{\a,j} \,\mbox{Im}H_{\a,j}
\eeq
is globally well-defined on $\M_{g,N}({\bf n})$. The critical points of this function on the leaves of the foliations described above are then described by the following theorem, which is the $N$-point generalization of Th. \ref{coords3}.
\begin{theo} Let $q_{\alpha}(t)$ be polynomials of degrees $m_{\alpha}$, such that $q_{\alpha}(t)$ does not have a constant term for $\alpha=2,\ldots,N$. Let $R_{\alpha}^Q$ be imaginary numbers adding up to zero. Let $(\Gamma,P_{\alpha},[z_{\alpha}])$ be a point in $\mathcal{M}_{g,1}({\bf n})$, which we assume to be embedded in $\mathcal{M}_{g,1}^m(n)$ via the data $q_{\alpha}(z)$ and $R_{\alpha}^Q$. Let $E$ and $Q$ be the real-normalized Abelian integrals defined by equations (\ref{defe3}), (\ref{defq31}) and (\ref{defq3a}), and let $D$ denote the intersection of the zero-divisors of $dE$ and $dQ$.

\noindent {\bf A.} The point $(\Gamma,P_{\alpha},[z_{\alpha}])$ is a critical point of $\H$ restricted to $\L$ if and only if

(i) There exists a meromorphic function $Y$ on $\Gamma$ with poles at $P_{\alpha}$ with principal parts
\beq\label{th321a}
Y=i\sum_{j=1}^{M_1} d_{1,j}z_1^{-j}+O(z_1)\mbox{ near }P_1,\quad
Y=i\sum_{j=1}^{M_{\alpha}} d_{\alpha,j}z_{\alpha}^{-j}+O(1)\mbox{ near }P_{\alpha},\quad \alpha=2,\ldots,N,
\eeq
and poles of order $\mu_s$ at the common zeros $Q_s$ of the differentials $dE$ and $dQ$
on $\Gamma$, where $\mu_s$ is the multiplicity of the common zero $Q_s$, and no other singularities.

\noindent {\bf B.} The point $(\Gamma,P_{\alpha},[z_{\alpha}])$ is a critical point of $\H$ restricted to $\L_E$ if and only if there exists a function satisfying condition (i) and in addition the following condition:

(ii) the differential $YdE$ is real-normalized, i.e.
\beq\label{th322}
\mbox{Im} \oint_{\g} YdE=0\quad\forall\gamma\in H_1(\Gamma,\mathbb{Z}).
\eeq

\noindent {\bf C.} The point $(\Gamma,P_{\alpha},[z_{\alpha}])$ is a critical point of $\H$ restricted to $\L_Q$ if and only if there exists a function satisfying condition (i) and in addition the following condition:

(iii) the differential $YdQ$ is real-normalized, i.e.
\beq\label{th323}
\mbox{Im} \oint_{\g} YdQ=0\quad\forall\gamma\in H_1(\Gamma,\mathbb{Z}).
\eeq

\noindent {\bf D.} The point $(\Gamma,P_{\alpha},[z_{\alpha}])$ is a critical point of $\H$ restricted to $\Lambda$ if and only if there exists a function satisfying conditions (i), (ii) and (iii).

\noindent {\bf E.} The point $(\Gamma,P_{\alpha},[z_{\alpha}])$ is a critical point of $\H$ restricted to $\Lambda'$ if and only if there exists a function satisfying condidions (i), (ii) and (iii), and such that near $P_{\alpha}$ we have
\beq\label{th321e}
Y=i\sum_{j=1}^{M_{\alpha}} d_{\alpha,j}z_{\alpha}^{-j}+O(z_{\alpha})\mbox{ near }P_{\alpha}.
\eeq
\label{main3}
\end{theo}

The proof of the theorem is similar to the one-point case. We introduce an extension of the homomorphism (\ref{dmap2}) to the $N$-point case:
\beq\label{dmap3}
f: T_P(\Lambda')\rightarrow \Omega(\Gamma,P_{\alpha},E,Q).
\eeq
For a tangent vector $\delta \in T_P(\Lambda')$, the differential $f(\delta)=\Omega_{\delta}=(\delta Q|_E) dE$ has simple poles at $P_1,\ldots,P_N$ and jumps along the cuts that are real linear combinations of $dE$ and $dQ$. The space of such differenials, which we denote by $\Omega(\Gamma,P_{\alpha},E,Q)$, has real dimension $6g+2N-2$. The image of the tangent subspace $T(\Lambda)$ is the subspace $\Omega(\Gamma,P_1,E,Q)$ of differentials having poles only at the point $P_1$.

The corresponding generalization of Prop. \ref{map2} is

\begin{prop}
Let $D=\sum\mu_sQ_s$ denote the intersection of the zero divisors of $dE$ and $dQ$. Then the image of the tangent space $T_P(\Lambda')$ under the homomorphism $f$ of (\ref{dmap3}) is the subspace $\Omega_D(\Gamma,P_{\alpha},E,Q) \subset\Omega(\Gamma,P_{\alpha},E,Q)$ of differentials having zeroes at $Q_s$ of orders $\mu_s$.  In particular, if the zero divisors of $dE$ and $dQ$ do not intersect, then the map (\ref{dmap3}) is an isomorphism.
\label{map3}
\end{prop}

The space $\Omega(\Gamma,P_{\alpha},E,Q)$ is spanned over $\mathbb{R}$ by the $6g$ differentials defined by conditions (\ref{normah})-(\ref{normq}), and in addition by the following differentials:

\begin{itemize} \item $2N-2$ single-valued meromorphic differentials $d\Omega_{\alpha}^R, d\Omega_{\alpha}^I$, $\alpha=2,\ldots,N$, which have real periods and have poles at $P_{\alpha}$ and $P_1$ with residues $\pm 1$ and $\pm i$, respectively:
\begin{equation}
\mbox{Res}_{P_\alpha} d\Omega_{\alpha}^R=1,\quad \mbox{Res}_{P_1} d\Omega_{\alpha}^R=-1,\quad
\mbox{Im}\oint_{\gamma}d\Omega_{\alpha}^R=0,\quad \gamma=A_i,B_i,
\end{equation}
\begin{equation}
\mbox{Res}_{P_\alpha} d\Omega_{\alpha}^I=i,\quad \mbox{Res}_{P_1} d\Omega_{\alpha}^I=-i,\quad
\mbox{Im}\oint_{\gamma}d\Omega_{\alpha}^I=0\quad \gamma=A_i,B_i.
\end{equation}
\end{itemize}

{\bf Proof of Theorem 3.2.}  Let $dY$ be the unique real-normalized Abelian differential with fixed branch at $P_1$ such that $Y$ has poles of the form (\ref{th321a}). Fix a real number $\rho>0$, and let $\gamma_{\alpha}$ denote a path of radius $\rho$ around the point $P_{\alpha}$ in the chosen local coordinate $z_{\alpha}$. Let $\delta\in T(\Lambda')$ be a tangent vector. A direct calculation then shows that
\begin{equation}
\delta H=\mbox{Re}\left(\displaystyle\sum_{\alpha=2}^N d_{\alpha,0}\delta T_{\alpha,0}\right)-
\displaystyle\sum_{\alpha=1}^N\mbox{Re}\,\mbox{Res}_{P_{\alpha}} Y \delta Q dE=
\mbox{Re}\left(\displaystyle\sum_{\alpha=1}^N d_{\alpha,0}\delta T_{\alpha,0}\right)-\frac{1}{2\pi}\mbox{Im}\displaystyle\oint_{\partial \Gamma}Y \Omega_{\delta},
\end{equation}
where $d_{\alpha,0}$ denotes the constant term of $Y$ at $P_{\alpha}$. For a tangent vector $\delta$ in $T(\Lambda)$ or in any smaller tangent subspace, we have $\delta T_{\alpha,0}=0$. Therefore, we can use the same reasoning as in Th. \ref{main2}, which proves parts A through D.

Finally, if $P$ is a critical point of $\mathcal{H}$ restricted $\Lambda'$, then $\delta \mathcal{H}=0$ for all variations of $\delta T_{0,\alpha}$, which can only happen if $d_{0,\alpha}=0$ for all $\alpha=2,\ldots,N$, which proves part E.


\begin{thebibliography}{99}
\bibitem{nov}
S.P. Novikov, ''A periodic problem for the Korteweg-de Vries equation'', {\it Funct. Anal. Appl.}, 8 (1974) n 3, 54-66

\bibitem{dmn}
B. Dubrovin, V. Matveev, S. Novikov, ``Non-linear equations of
Korteweg-de Vries type, finite zone linear operators and Abelian varieties'',
Uspekhi Mat. Nauk {\bf 31}:1 (1976) 55-136.

\bibitem{non-statShrod}
I.M. Krichever, {\it The spectral theory of ``finite-gap'' nonstationary Schrödinger operators},
The nonstationary Peierls model. (Russian) Functional. Anal. Appl. 20 (1986), no. 3, 42--54, 96.

\bibitem{boutroux}
Boutroux P. {\it Recherches sur les transcendentes de M Painlev\'e et l'edute asymtotique des \'equations diff\'erentielles du second orddre}, Ann Sci Ecol.Norm Sup\'er 30 (1913) 255-376, 31(1914), 99-159

\bibitem{kr-pen}
Krichever I.M. {\it On Heizenberg relations for the ordinary linear differential operators},
IHES preprint 1990, Bur-Sur-Yvette

\bibitem{gamba}
Fucito F., Gamba A., Martellini M., Ragnisco O., {\it Non-linear WKB analysis of the string equation},  
International Journal of Modern Physics B, 6(11-12), pp. 2123--2147, 1992.

\bibitem{grinevich-novikov}
Grinevich P. G., Novikov S. P., ``String equation---2. Physical solution'', {\it St. Petersburg Math. J.} 6, p. 553, 1995.

\bibitem{kapaev}
Kapaev A. A., ``Monodromy approach to the scaling limits in the isomonodromy systems'', arXiv.org:nlin/0211022

\bibitem{david1}
David F., ``Phases of the large-N matrix model and non-perturbative effects in 2D gravity'', {\it Nuclear Phys. B} 348(4), pp. 507--524, 1991.

\bibitem{david2}
David F., ``Non-perturbative effects in matrix models and vacua of two dimensional gravity'', {\it Phys. Lett. B} 302(4), pp. 403--410, 1993.

\bibitem{bertola-mo}
Bertola M, Mo M. Y., ``Commuting difference operators, spinor bundles and the asymptotics of orthogonal polynomials with respect to varying complex weights'', arXiv.org:math-ph/0605043, 2006.

\bibitem{wiegmann}
Lee, S.Y., Teodorescu, R., Wiegmann, P., ``Shocks and finite-time singularities in Hele-Shaw flow'', {\it Physica D: Nonlinear Phenomena},  238(14), p. 1113-1128.



\bibitem{kr-tau} I. Krichever: {\it The $\tau$-function of the universal Whitham hierarchy, matrix models,
and topological field theories}, Comm. Pure Appl. Math. {\bf 47} (1994), 437--475.

\bibitem{kp1} I. Krichever, D.H. Phong: {\it On the integrable geometry of $N=2$ supersymmetric gauge
theories and soliton equations}, J. Differential Geometry {\bf 45} (1997) 445-485.

\bibitem{kp2} I. Krichever, D.H. Phong: {\it Symplectic forms in the theory of solitons},
Surveys in Differential Geometry {\bf IV} (1998), edited by C.L. Terng and K. Uhlenbeck, 239-313,
International Press.

\bibitem{kr-grush1}
S. Grushevsky, I. Krichever, {\it The universal Whitham hierarchy and the geometry of
the moduli space of pointed Riemann surfaces}, arXiv:0810.2139


\bibitem{kr-grush2}
S. Grushevsky, I. Krichever, {\it Vanishing of tautological classes on the moduli
space of curves}, in preparation



\bibitem{dzhamay}
Dzhamay, A., ''Real-normalized Whitham hierarchies and the WDVV equations'', Internat. Math. Res. Notices, 2000, no. 21, 1103--1130.





\end{thebibliography}
\end{document}